\documentclass[11pt,twoside]{article}

\usepackage{amsbsy,amsfonts,amsmath,amssymb,eucal,mathrsfs}
\usepackage[all]{xy}
\usepackage{pstricks}
\usepackage{pst-plot}

\newtheorem{definition}{Definition}[section]

\newtheorem{prop}[definition]{Proposition}

\newtheorem{lemm}[definition]{Lemma}

\newtheorem{coro}[definition]{Corollary}
\newtheorem{theo}[definition]{Theorem}
\newtheorem{notation}[definition]{Notation}

\newtheorem{construction}[definition]{Construction}

\newtheorem{remark}[definition]{Remark}

\newtheorem{remarks}[definition]{Remarks}

\newtheorem{example}[definition]{Example}

\newtheorem{examples}[definition]{Examples}

\newtheorem{nothing}[definition]{$\!\!$}

\newtheorem{definition*}{Definition}[section]
\newenvironment{defi*}{\begin{definition*} \rm}{\end{definition*}}
\newtheorem{definitions*}[definition*]{Definitions}
\newenvironment{defis*}{\begin{definitions*} \rm}{\end{definitions*}}
\newtheorem{prop*}[definition*]{Proposition}
\newtheorem{lemm*}[definition*]{Lemma}
\newtheorem{coro*}[definition*]{Corollary}
\newtheorem{theo*}[definition*]{Theorem}
\newtheorem{remark*}[definition*]{Remark}
\newenvironment{rema*}{\begin{remark*} \rm}{\end{remark*}}
\newtheorem{remarks*}[definition*]{Remarks}
\newenvironment{remas*}{\begin{remarks*} \rm}{\end{remarks*}}
\newtheorem{example*}[definition*]{Example}
\newenvironment{exam*}{\begin{example*} \rm}{\end{example*}}
\newtheorem{examples*}[definition*]{Examples}
\newenvironment{exams*}{\begin{examples*} \rm}{\end{examples*}}
\newtheorem{nothing*}[definition*]{$\!\!$}
\newenvironment{noth*}{\begin{nothing*} \rm}{\end{nothing*}}

\begin{document}

\def\pt{\{{\rm pt}\}}
\def\ra{\rightarrow}
\def\s{\sigma}\def\OO{\mathbb O}\def\PP{\mathbb P}\def\QQ{\mathbb Q}
 \def\CC{\mathbb C} \def\ZZ{\mathbb Z}\def\JO{{\mathcal J}_3(\OO)}
\newcommand{\G}{\mathbb{G}}
\def\proof{\noindent {\it Proof.}\;}
\def\qed{\hfill $\square$} 
\def \uh {{\widehat{u}}}
\def \vh {{\widehat{v}}}
\def \wh {{\widehat{w}}}
\def \xh {{\widehat{x}}}
\def \cS {\mathcal S}
\def \ft {\mathfrak t}
\def \fp {\mathfrak p}
\def \fe {\mathfrak e}
\def \ff {\mathfrak f}
\def \fsl {\mathfrak {sl}}
\def \fso {\mathfrak {so}}
\def \fspin {\mathfrak {spin}}
\def \fg {\mathfrak g}
\def\s{\sigma}
\def\trace{\mathrm{trace}}



\newcommand{\N}{\mathbb{N}}
\newcommand{\Z}{\mathbb Z}
\newcommand{\R}{\mathbb{R}}
\newcommand{\Q}{\mathbb{Q}}
\newcommand{\C}{\mathbb{C}}
\renewcommand{\H}{\mathbb{H}}
\renewcommand{\O}{\mathbb{O}}
\newcommand{\F}{\mathbb{F}}
\newcommand{\norme}[1]{|| #1 ||}
\newcommand{\scal}[1]{\langle #1 \rangle}

\renewcommand{\a}{{\cal A}}
\newcommand{\az}{\a_\Z}
\newcommand{\ak}{\a_k}

\newcommand{\rc}{\R_\C}
\newcommand{\cc}{\C_\C}
\newcommand{\hc}{\H_\C}
\newcommand{\oc}{\O_\C}

\newcommand{\rk}{\R_k}
\newcommand{\ck}{\C_k}
\newcommand{\hk}{\H_k}
\newcommand{\ok}{\O_k}

\newcommand{\rz}{\R_Z}
\newcommand{\cz}{\C_Z}
\newcommand{\hz}{\H_Z}
\newcommand{\oz}{\O_Z}

\newcommand{\RR}{\R_R}
\newcommand{\CR}{\C_R}
\newcommand{\HR}{\H_R}
\newcommand{\OR}{\O_R}

\newcommand{\re}{\mathtt{Re}}


\newcommand{\dual}{{\bf v}}
\newcommand{\com}{\mathtt{Com}}
\newcommand{\rg}{\mathtt{rg}}

\newcommand{\g}{\mathfrak g}
\newcommand{\h}{\mathfrak h}
\renewcommand{\u}{\mathfrak u}
\newcommand{\n}{\mathfrak n}
\newcommand{\e}{\mathfrak e}
\newcommand{\plie}{\mathfrak p}
\newcommand{\q}{\mathfrak q}
\newcommand{\liesl}{\mathfrak {sl}}
\newcommand{\so}{\mathfrak {so}}
\newcommand{\p}{\PP}

 \title{Quantum cohomology of minuscule homogeneous spaces
   III\linebreak Semi-simplicity and consequences}
 \author{P.E. Chaput, L. Manivel, N. Perrin}

\maketitle

\begin{abstract}
We prove that the quantum cohomology ring of any minuscule or
cominuscule homogeneous space, specialized at $q=1$, is semisimple. 
This implies that complex conjugation defines an algebra automorphism
of the quantum cohomology ring localized at the quantum
parameter. We check that this involution coincides with the strange 
duality defined in \cite{cmp2}. We deduce Vafa-Intriligator type 
formulas for the Gromov-Witten invariants. 
\end{abstract}

 {\def\thefootnote{\relax}
 \footnote{ \hspace{-6.8mm}
 Key words: quantum cohomology, minuscule homogeneous spaces,
  Schubert calculus, quantum Euler class. \\
 Mathematics Subject Classification: 14M15, 14N35}
 }

\section{Introduction}

This paper is a sequel to \cite{cmp} and \cite{cmp2}, where we develop a 
unified approach to the quantum cohomology of (co)minuscule homogeneous
manifolds $X=G/P$. 

Recall that a $\ZZ$-basis for the ordinary
cohomology ring $H^*(X)$ (or for the Chow ring $A^*(X)$) of $X$ 
is given by the Schubert classes
$\sigma(w)$, where $w\in W_X$ belongs to the set of minimal lengths
representatives of $W/W_P$, the quotient of the Weyl group $W$ of $G$
by the Weyl group $W_P$ of $P$. 
The Schubert classes are also a basis over $\ZZ[q]$ of the (small) 
quantum Chow ring $QA^*(X)$, whose associative product is
defined in terms of 3-points Gromov-Witten invariants. If $R$ is a ring,
we denote by 
$QA^*(X,R)$ the tensor product $QA^*(X) \otimes_\ZZ R$ and
$QA^*(X,R)_{loc}$ its localization at $q$, that is, 
$QA^*(X,R)_{loc}=QA^*(X,R)\otimes_{\ZZ[q]}\ZZ[q,q^{-1}].$
The main result of \cite{cmp2}, {\it strange duality}, 
was that one could define, 
for any $w\in W_X$, a non negative integer $\delta(w)$, and an algebraic number
$\zeta(w)$, in such a way that the endomorphism
$\iota$ of $QA^*(X,\R)_{loc}$, defined  by 
$$\iota(q)=q^{-1} \quad \mathrm{and}  \quad 
\iota(\sigma(w))=q^{-\delta(w)}\zeta(w)\sigma(\iota(w)),$$
be a ring involution. 

In this paper we give a natural explanation of the existence of 
such an involution, which turns out to be directly related 
to the semi-simplicity 
of the quantum cohomology ring specialized at $q=1$. For the classical cases 
this semisimplicity has been known for some time, as it can readily be read
off the explicit presentations that have been found. We complete the picture 
by checking the semisimplicity in the exceptional cases. 
  
To this end, we recall that for $X=G/P$ a minuscule homogeneous space, the 
rational quantum cohomology ring can be described as
$$QH^*(G/P)=\QQ[\ft]^{W_P}[q]/(I_{d_1},\ldots ,I_{d_m}-q),$$
where $\ft$ is a Cartan subalgebra of the Lie algebra $\fg$ of $G$, and $I_{d_1},\ldots ,I_{d_m}$
are homogeneous generators of the $W$-invariants;
the maximal degree  $d_m$  is the 
Coxeter number of $\fg$. 

Let $Z(\fg)$ denote the subscheme of $\ft$ defined by the equations $I_{d_1}=\cdots =I_{d_{m-1}}=I_{d_m}-1=0$.
Observe that it does not depend (up to homothety) on the choice of the invariants. Moreover, it is
a finite scheme of length $\# W$, since if we replace $I_{d_m}-1$ by $I_{d_m}$ in this set of equations, 
we get the spectrum of the Chow ring of the full flag variety of $G$, which is a finite scheme of length $\# W$
supported at the origin. The following result indicates a fundamental difference between the classical and the 
quantum settings.

\begin{prop}\label{free}
For any simple Lie algebra $\fg$ (except possibly $\ff_4$ and $\fe_8$),  
the scheme $Z(\fg)$ is reduced, and is a free $W$-orbit. 
\end{prop}

In fact this result is relevant for quantum cohomology only in type $A,D,E$, with $E_8$
excepted. For the classical types it is in fact very easy to check, but the cases of $E_6$ 
and $E_7$ are somewhat more involved. It implies the ``minuscule'' part of our next statement:
 
\begin{coro}\label{reduced}
For any minuscule or cominuscule rational homogeneous space $G/P$, the 
algebra $QA^*(G/P)_{q=1}\otimes_\QQ\CC\simeq \CC[W/W_P]$ is semi-simple.
\end{coro}
 
More intrisically, the spectrum $Z(G/P)$ of the quantum algebra 
$QA^*(G/P)_{q=1}$ is $Z(\fg)/W_P$, at least in the minuscule case. 

Now, any commutative semi-simple finite-dimensional algebra $H$ is a product
of fields, hence over $\R$ it decomposes as $H = \R^n \oplus \C^p$. By conjugating 
the complex factors we get a canonical algebra automorphism, the complex 
involution.

We therefore get an algebra 
involution of $QA^*(G/P,\R)_{q=1}$. We point out that, 
{\it because $Z(G/P)$ is reduced}, 
this involution lifts to an algebra involution of $QA^*(G/P)_{loc}$, 
mapping $q$ to 
$q^{-1}$ and any class of degree $d$ to a class of degree $-d$
(see theorem \ref{complexinv}). 
This leads to a new interpretation of  strange duality. 

\begin{theo}\label{sameinv}
For any minuscule or cominuscule homogeneous space $X=G/P$, complex conjugation and 
strange duality define the same involution.
\end{theo}

The proof is case by case. In fact this result has already been observed by Hengelbrock 
for the case of Grassmannians, with a different method \cite{heng}. 
What we will check is that strange duality and the complex involution 
coincide on a set of generators of the quantum cohomology ring. In the classical cases we will also
provide a direct check that the complex involution is given by the same expressions as in \cite{cmp2}. 
In particular this will explain the occurrence of the  irrationalities introduced by the function $\zeta$  
(which one can eventually get rid of by rescaling $q$). 

The advantage of this approach of strange duality through complex conjugation is that, beside
of being conceptually enlightening, the fact that it is an algebra automorphism becomes 
completely obvious -- while in \cite{cmp2} this was the result of painful computations, 
especially in the exceptional cases. What is not clear a priori is that the complex 
conjugate of a Schubert class is again (a multiple of) a Schubert class, while in \cite{cmp2}
this was given by the very definition of strange duality. It would be interesting to have 
a conceptual explanation of that phenomenon (which does no longer hold true on non 
minuscule or cominuscule spaces).   

We stress that the smoothness of the finite scheme $Z(G/P)$ plays an essential role here. 
In section 7, we consider the case of $G_{\omega}(2,6)$, the Grassmannian of isotropic planes 
in a six-dimensional symplectic complex vector space. This is the simplest example of a 
homogeneous space with Picard number one which is neither minuscule nor cominuscule. 
We check that its quantum cohomology ring is not semisimple; in fact the scheme $Z(G_{\omega}(2,6))$
is made of ten simple points and one double point.  Moreover, the existence of this double 
point prevents the complex conjugation from being lifted to an involution of the quantum
cohomology ring reversing degrees. 

Finally, we use our schemes $Z(G/P)$ to obtain Vafa-Intriligator type formulas for the 
Gromov-Witten invariants. We express these formulas in a uniform way, in terms of a 
{\it quantum Euler class} $e(X)$ introduced in \cite{ab} for any projective manifold $X$. 
Abrams proved that the invertibility  of that class is equivalent to the semi-simplicity 
of the quantum cohomology ring $QH^*(X)$ (after specialization of the quantum parameters)
In the (co)minuscule setting, we prove that the quantum Euler class is simply given by the 
product of the positive roots of $\fg$ that
are not roots of $\fp=Lie(P)$.  The Vafa-Intriligator type formulas that we obtain
are expressed in terms of that class. They are equivalent to the formulas obtain in 
\cite{st} and \cite{rie} for Grassmannians, but they are simpler than the formulas 
given in \cite{cheong} for the other classical cases. 

\vskip 0.5 cm

{\bf Aknowledgement:} we would like to thank Konstanze Rietsch for
usefull discussions on several aspects of the quantum cohomology in
the classical cases.

\section{The complex involution}

Let $X$ be a projective variety with Picard number one. Let  
$QH^*(X,\C) = H^*(X,\C) \otimes_\C \C[q]$
be its small quantum cohomology ring over the
complex numberss, and
$QH^*(X,\C)_{loc}$ the algebra obtained by inverting the quantum parameter $q$.

\begin{theo}\label{complexinv}
Suppose that the spectrum of the finite dimensional algebra $QH^*(X)_{q=1}$ 
is a reduced finite scheme.   
Then there exists an algebra automorphism of $QH^*(X,\C)_{loc}$ mapping
$q$ to $q^{-1}$.
\end{theo}

\proof 
The inclusion $\C[q,q^{-1}] \hookrightarrow QH^*(X,\C)_{loc}$ 
yields a finite morphism
$\pi : C \rightarrow D$ of curves over $\C$. We consider the involution $i$ of
$D$ given by $q \mapsto q^{-1}$; the theorem states that we can lift
$i$ to $C$, under the hypothesis that $C$ be a smooth curve. 

Let us consider a homogeneous presentation 
$H^*(X,\R) = \R[X_1,\ldots,X_n]/(r_1,\ldots,r_k)$ 
of $H^*(X,\R)$. The quantum cohomology ring can be presented as
$\R[X_1,\ldots,X_n,q]/(R_1,\ldots,R_k)$, where
$R_i$ is again a homogeneous relation, that specialises to $r_i$ when $q=0$.

If $(x_i,q) \in C(\C)$ is a complex point of $C$, 
then $(\overline x_i,\overline q) \in C(\C)$
because $C$ is defined over $\R$, and thus, by homogeneity, 
$(\overline x_i/\norme q^{2\deg(q)/\deg(x_i)},\overline q/\norme q^2)$
also belongs
to $C(\C)$. This is a complex  point of $C$ over
$\overline q/\norme q^2 = q^{-1}$. 
We claim that the map
$$j : (x_i,q) \mapsto (\overline x_i/\norme q^{2\deg(q)/\deg(x_i)},
\overline q/\norme q^2)$$
is algebraic. Indeed, consider the fiber product
$C \times_D C$, where the first morphism $C \rightarrow D$ 
is $\pi$, and the second
$i \circ \pi$. Let $C_0$ denote the connected component in
$(C \times_D C)(\C)$, given as the set of pairs $((x_i,q),j(x_i,q))$. It is
algebraic, and the morphism $C_0 \rightarrow C$ induced by the first projection
is finite of
degree one, thus an isomorphism. So the theorem is proved.
\qed

\medskip
We call this involution $i$ of $QH^*(X,\C)_{loc}$ 
the {\it complex involution}.

\section{Grassmannians}

Let $G(d,n)$ denote the Grassmannian of $d$-dimensional subspaces of an $n$-dimensional
vector space. Its quantum cohomology ring can be described as 
$$QA^*(G(d,n))=\ZZ[x_1,\ldots ,x_n]^{\cS_d\times\cS_{n-d}}[q]/(e_1,\ldots ,e_{n-1}, e_n-q),$$
where $e_1, \ldots ,e_n$ are the elementary symmetric functions in the $n$ indeterminates
$x_1,\ldots ,x_n$ \cite{st}. Here the symmetric groups $\cS_d$ and  $\cS_{n-d}$ act by 
permutation of the first $d$ and last $n-d$ variables, so we only consider symmetric functions 
in these two set of variables. Usually, the relations $e_1,\ldots ,e_{n-1}, e_n$ are used to 
eliminate one of these two sets of variables, but we will not do that.

\medskip\noindent {\it Proof of \ref{free}}. 
The equations defining $Z(\fsl_n)$ in $\CC^n$ 
are $e_1=\cdots = e_{n-1}=e_n-1=0$.  Thus $Z(\fsl_n)$ is 
the set of $n$-tuples $(\zeta_1,\ldots ,\zeta_n)$ of distinct $n$-th roots of 
$(-1)^{n-1}$. This is certainly a 
free orbit of the symmetric group $\cS_n$. \qed

\medskip
We can therefore interpret the quantum cohomology ring of $G(d,n)$ at $q=1$, 
as 
$$QH^*(G(d,n))_{q=1}=\QQ[Z(G(d,n))],$$
where $Z(G(d,n))$ denotes the set of (unordered) $d$-tuples of distinct $n$-th roots 
of $(-1)^{n-1}$.

\medskip\noindent {\it Proof of \ref{sameinv}}. 
As an algebra, $QH^*(G(d,n))_{q=1}$ is generated by the special Schubert classes $\s(k)$, 
whose corresponding functions on $Z(G(d,n))$ are the $k$-th elementary symmetric functions. 
What we need to prove is that the complex involution maps such a special Schubert class $\s(k)$ 
to the class $\sigma(n-d,1^{d-k}).$ To check this we simply compute the complex conjugate 
of $\s(k)$ as follows:
\begin{eqnarray*}
e_k(\bar\zeta_1,\ldots ,\bar\zeta_d) &= &e_k(\zeta^{-1}_1,\ldots ,\zeta^{-1}_d) \\
          &= &(\zeta_1\ldots \zeta_d)^{-1}e_{d-k}(\zeta_{1},\ldots ,\zeta_d) \\
          &= &(\zeta_{d+1}\ldots \zeta_n)e_{d-k}(\zeta_{1},\ldots ,\zeta_d).
\end{eqnarray*}
Observe that the function $\zeta_{d+1}\ldots \zeta_n=e_{n-d}(\zeta_{d+1},\ldots ,\zeta_n)=
h_{n-d}(\zeta_{d+1},\ldots ,\zeta_n)$. Therefore $\bar e_k=h_{n-d}e_{d-k}=s_{n-d,1^{d-k}}$. \qed  

\medskip
For completeness we deduce the complex conjugate of any Schubert class and recover the 
formulas given by Postnikov \cite{post}. 

\begin{prop} The complex conjugate of the Schubert class $\sigma (\lambda)$, is the 
 Schubert class $\sigma (\iota(\lambda))$. 
\end{prop}

Here $\iota(\lambda)$ denotes the partition deduced  from $\lambda$ by a 
simple combinatorial process (see \cite{post} and \cite{cmp2}). Recall that $\lambda$ is a partition 
whose diagram is contained in a $d\times (n-d)$ rectangle. Let $c$ be the size of the 
Durfee square of $\lambda$, that is, the largest integer such that $\lambda_c\ge c$. 
Write $\lambda$ as $(c+\mu,\nu)$, where now $\mu$ is contained in a $c\times (n-c)$ rectangle
and $\mu$ in a $(d-c)\times c$ rectangle. Denote by $p(\mu)$ and $p(\nu)$ the complementary
partitions in these respective rectangles. Then $\iota(\lambda)=(c+p(\mu),p(\nu))$. 

\medskip\proof 
We use the the fact that the Giambelli formulas
hold in the quantum cohomology ring, as proved by Bertram in  \cite{bert}. Thus, for any partition $\lambda$, 
\begin{eqnarray*}
\bar\sigma (\lambda) &= &(h_{n-d})^{n-d}\det (\sigma(d-\lambda^*_i+i-j))_{1\le i,j\le n-d}\\
          &= &(h_{n-d})^{n-d}\det (\sigma(d-\lambda^*_{n-d+1-i}-i+j))_{1\le i,j\le n-d}\\
          &= &(h_{n-d})^{n-d}\sigma (p(\lambda)),
\end{eqnarray*}
where $\sigma (p(\lambda))$ is the Schubert class Poincar{\'e} dual to $\sigma (\lambda)$. 
But $h_{n-d}$ is invertible in the quantum cohomology ring, with inverse $\sigma(d)$, 
and  $\sigma(d)^{n-d}$ is the punctual class $\s(pt)$. Since the multiplication by the punctual
class is given (for $q=1$) by the formula $\s(pt)*\sigma(\mu)=\sigma(p\iota(\mu))$ (see \cite{cmp2}, 
Theorem 3.3), we finally get 
$$\bar\sigma (\lambda)=\s(pt)^{-1}*\sigma (p(\lambda))=\sigma (\iota(\lambda)).$$
Note the interesting fact that we obtain directly the relation that holds for any (co)minuscule 
homogeneous space, between the complex conjugation, Poincar{\'e} duality, and the quantum product 
by the punctual class. \qed 

\section{Orthogonal Grassmannians and Quadrics}

\subsection{Orthogonal Grassmannians}
Let $G_Q(n+1,2n+2)$ denote the orthogonal Grassmannian, that is, one of the two families 
of maximal isotropic subspaces in some vector space of dimension $2n+2$ endowed with a 
non degenerate quadratic form. Its quantum cohomology ring can be described as 
$$QA^*(G_Q(n+1,2n+2))=\ZZ[x_1,\ldots ,x_{n+1}]^{\cS_{n+1}}[q]/(E_1,\ldots ,E_{n-1}, E_n-2q,e_{n+1}),$$
where $E_1, \ldots ,E_n,E_{n+1}=e_{n+1}^2$ are now the elementary symmetric functions in the squares of 
the $n+1$ indeterminates $x_1,\ldots ,x_{n+1}$ (\cite{BKT2}, Theorem 1).

\medskip\noindent {\it Proof of \ref{free}}. 
The equations defining $Z(\fso_{2n+2})$ 
are $E_1=\cdots = E_{n-1}=E_{n}-4=e_{n+1}=0$. 
The set of solutions of this equation is  
$$Z(\fso_{2n+2})=\{(\zeta_1,\ldots ,\zeta_k,0,\zeta_{k+1},\ldots ,\zeta_n)\},$$
where the squares of $\zeta_1,\ldots ,\zeta_n$ are the $n$ distinct $n$-th roots of $(-1)^{n-1}.4$.
This is certainly a free orbit of the Weyl group $W(D_{n+1})=\cS_{n+1}\times \ZZ_2^n$. \qed

\medskip
We can therefore interpret the quantum cohomology ring of $G_Q(n+1,2n+2)$ at $q=1$, as 
$$QH^*(G_Q(n+1,2n+2))_{q=1}=\QQ[Z(G_Q(n+1,2n+2))],$$
where $Z(G_Q(n+1,2n+2))$ is identified with
the set of (unordered) $n$-tuples of square roots of 
the $n$ distinct $n$-th roots of $(-1)^{n-1}.4$.
Note that $\# Z(G_Q(n+1,2n+2))=2^n$, as expected. 

\medskip\noindent {\it Proof of \ref{sameinv}}. 
The quantum cohomology ring is generated by the special Schubert classes $\s(k)$,  
for $1\le k\le n$. The class $\s(k)$ is represented by 
$e_k/2$, where $e_k$ is the $k$-th elementary symmetric function
in $x_1,\ldots ,x_n$ (see \cite{BKT2}, Theorem 1).  What we need to prove is that the complex involution maps 
a special Schubert class $\s(k)$ to a suitable multiple of the class $\sigma(n,n-k).$
We compute the complex conjugate of $\s(k)$ as follows:
\begin{eqnarray*}
e_k(\bar\zeta_1,\ldots ,\bar\zeta_n) &= &4^{k/n}e_k(\zeta^{-1}_1,\ldots ,\zeta^{-1}_n) \\
          &= &(\zeta_1\ldots \zeta_n)^{-1}e_{n-k}(\zeta_{1},\ldots ,\zeta_n) \\
          &= &4^{k/n-1}(\zeta_{d+1}\ldots \zeta_n)e_{n-k}(\zeta_{1},\ldots ,\zeta_n),
\end{eqnarray*}
where we have used the fact that $(\zeta_{d+1}\ldots \zeta_n)^2=e_{n}(\zeta_{1},\ldots ,\zeta_n)^2=4$. 
Otherwise stated, $$\bar \sigma(k)=2^{2k/n-1}\sigma(n)\sigma(n-k)=2^{2k/n-1}\sigma(n,n-k),$$
in complete agreement with Proposition 4.7 in \cite{cmp2}. 
Note finally that since $e_n^2=4$, $e_n$ is real, hence $\sigma(n)$ is fixed by the complex involution, 
and we are  done. \qed

\medskip
We deduce that the complex conjugate of any Schubert class $\s(\lambda)$ is given by a suitable 
multiple of $\s(\iota(\lambda))$, where $\iota(\lambda)$ is defined as follows. Write $\lambda=
(\lambda_1>\cdots >\lambda_{2\delta(\lambda)})$, ending with a zero part if necessary. Then 
$\iota(\lambda)=(n-\lambda_{2\delta(\lambda)}>\cdots >n-\lambda_1)$. Let  $z(\lambda)=\frac{2|\lambda|}{n}
-(\ell(\lambda)+\delta_{\lambda_1,n})$. 

\begin{prop}
The complex involution maps the Schubert class $\s(\lambda)$ 
to $2^{z(\lambda)}\sigma(\iota(\lambda))$.
\end{prop}

\proof
We first consider classes $\sigma(i,j)$, with 
$i>j>0$. Suppose for example that $i+j>n$; then
$$\sigma(i,j)=\sigma(i)\sigma(j)+2\sum_{k=1}^{n-i-1}(-1)^k\sigma(i+k)\sigma(j-k)+
(-1)^{n-i}\sigma(n)\sigma(i+j-n).$$
Using the fact that $\sigma(n)^2=1$, we deduce that 
$$\bar\sigma(i,j)=2^{\frac{2i+2j}{n}-2}\{\sigma(n-j)\sigma(n-i)+2\sum_{k=1}^{n-i-1}(-1)^k\sigma(n-j+k)
\sigma(n-i-k)+(-1)^{n-i}\sigma(2n-i-j)\},$$
that is, $\bar\sigma(i,j)=2^{\frac{2i+2j}{n}-2}\sigma(n-j,n-i)$. A similar computation 
shows that this formula also holds for $i+j\le n$. 
But then the Giambelli type formula (see \cite{BKT1}, Theorem 7)
$$\sigma(\lambda)=\mathrm{Pfaff}(\sigma(\lambda_i,\lambda_j))_{1\le i<j\le 2\delta(\lambda)}$$
immediately implies that $\bar\sigma(\lambda)= 2^{z(\lambda)}\sigma(\iota(\lambda))$. \qed

\subsection{Quadrics}

We will only treat the case of quadrics of even dimension (which are minuscule). The case
of quadrics of odd dimensions (which are cominuscule) is very similar. 

So let $\QQ^{2n}$ be a quadric of even dimension $2n$. There are two Schubert classes 
$\s_+$ and $\s_-$ in middle dimension $n$, and a single one $\s_k$ in every other dimension $k\ne n$. 
In terms of the hyperplane class $H=\s_1$, one has $\s_k=H^k$ for $k\le n-1$, $\s_k=H^k/2$ for $k\ge n+1$,
and $H^m=\s_++\s_-$ (in the classical Chow ring). See \cite{cmp}, section 4.1, for more details. 

\smallskip
We have just seen that $Z(\fso_{2n+2})$ is the set of $(n+1)$-tuples 
$(\zeta_1,\ldots ,\zeta_k,0,\zeta_{k+1},\ldots ,\zeta_n),$
where the squares of $\zeta_1,\ldots ,\zeta_n$ are the $n$ distinct $n$-th roots of $4$. 
The Weyl group $W_P$ of the parabolic $P$ defining $\QQ^{2n}$ is the fixator of the first coordinate. 
It has $2n+2$ orbits in $Z(\fso_{2n+2})$. First, there are $2n$ orbits $O(\zeta)$ defined by their 
non zero first coordinate $\zeta$, which can be any $2n$-th root of $4$. Second, there are two orbits 
$O(+)$ and $O(-)$ with zero first coordinate, and defined by the fact that the product of the non 
zero coordinates is $\pm 2$.
Since $2n+2$ is also the dimension of $H^*(\QQ^{2n})$, this confirms that $Z(\QQ^{2n})$ is reduced 
and $QA^*(\QQ^{2n})_{q=1}$ is semi-simple. 

The $W_P$-invariant polynomials are generated by $t_0$, the first coordinate, and the product 
$P=t_1\cdots t_n$ of the other coordinates.  
The algebra $QA^*(\QQ^{2n})_{q=1}$, considered as an algebra of
functions on the set $Z(\QQ^{2n})$, is determined by the following table:
$$\begin{array}{lccc}
 & H & H^n & P \\
 & & & \\
O(\zeta)\hspace*{3mm} & \zeta & \zeta^n& 0 \\
O(+) & 0 & 0  & 2 \\
O(-) & 0 & 0  & -2 
\end{array}$$
Observe that the two degree $n$ classes $H^n$ and $P$ are real. Moreover we easily get that 
$\overline H=4^{\frac{1-n}{n}}H^{2n-1}$. There remains to express $\s_+$ and $\s_-$ in terms
of $H^n$ and $P$. Using the formulas of \cite{cmp2}, we see that this expression depends 
on the parity of $n$; we get 
$$\s_{\pm}=\frac{1}{2}(H^n\pm iP)\quad \mathrm{for}\; n\; \mathrm{even}, \qquad
\s_{\pm}=\frac{1}{2}(H^n\pm P)\quad \mathrm{for}\; n\; \mathrm{odd}.$$
Thus $\overline\sigma_\pm=\sigma_\mp$ for $n$ even, but $\overline\sigma_\pm=\sigma_\pm$ for 
$n$ odd. Comparing with \cite{cmp2}, we see that strange duality coincides with the complex involution. 
 

\section{Lagrangian Grassmannians}
Let $G_{\omega}(n,2n)$ denote the Lagrangian Grassmannian parametrizing maximal isotropic subspaces 
in some symplectic vector space of dimension $2n$. This is a cominuscule, but not a minuscule homogeneous space. 
Its quantum cohomology ring can be described as 
$$QA^*(G_{\omega}(n,2n))=\ZZ[x_1,\ldots ,x_{n}]^{\cS_{n}}[q]/(R_1,\ldots ,R_{n}),$$
where the relations are $R_k=E_k-(-1)^{2k-n-1}qe_{2k-n-1}$ for $1\le k\le n$ \cite{T}. 
Recall that $e_k$ (resp. $E_k$) is the $k$-th elementary symmetric functions of the 
variables $x_1,\ldots ,x_{n}$ (resp. of their squares). Moreover, we used the convention that 
$e_k=0$ for negative $k$, so that the relations have no quantum correction in degree less than or equal 
to $n+1$. All the higher degree relations involve a quantum term, in contrast with the minuscule 
case for which only the relation of highest degree receives a quantum correction.  This makes that 
the quantum relations are no longer $W$-invariants, so 
there is no point in proving \ref{free}. 
Instead, we directly prove that $Z(G_{\omega}(n,2n))$ is reduced, hence the semi-simplicity 
of the quantum cohomology ring at $q=1$. 

\medskip\noindent {\it Proof of \ref{reduced}}. 
Observe that the relations take a simpler form if we introduce formally an auxiliary variable 
$x_0$ and let $q=e_{n+1}$. The point is that the $R_k$'s are then formally the same as the 
$k$-th elementary symmetric functions of the squares of $x_0,x_1,\ldots ,x_{n}$. This implies that 
the spectrum of the quantum cohomology ring at $q=1$ is supported on the set of (unordered) 
$(n+1)$-tuples $(\zeta_0,\zeta_1,\ldots,\zeta_n)$ such that $\zeta_0\zeta_1\cdots\zeta_n=1$
and $\zeta_0^2,\zeta_1^2,\ldots ,\zeta_n^2$ are the $(n+1)$-th different roots of $(-1)^n$. 
There are exactly $2^n$ such unordered $(n+1)$-tuples, as was already observed in \cite{cheong}. 
Since this coincides with the 
dimension of $QA^*(G_{\omega}(n,2n))_{q=1}$ as a vector space, this algebra is semi-simple. \qed

\medskip\noindent {\it Proof of \ref{sameinv}}. 
We leave this to the reader, since this case is even easier than the previous ones. Indeed the complex 
conjugation sends a special Schubert class $\s(k)$ to a suitable multiple of another special 
Schubert class $\s(n+1-k)$. A consequence is that the fact that this defines an involution of the 
quantum cohomology can be directly checked on the relations: $R_k$ is simply changed into a 
multiple of $R_{n+1-k}$, and the claim follows. \qed

\section{The exceptional cases}

\subsection{The Cayley plane}

We first recall some facts on the quantum cohomology ring of the Cayley plane $\OO\PP^2=E_6/P_1$. 
This is a sixteen dimensional variety, of index $12$. The Schubert classes are organized in 
the following Hasse diagram:

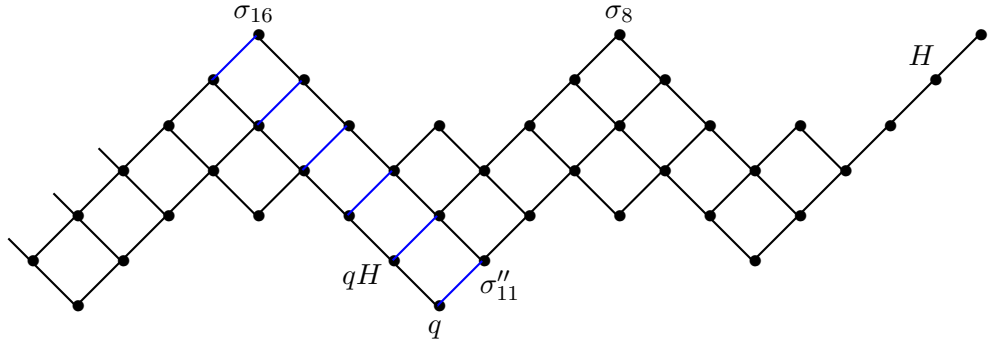
\begin{center}
\psset{unit=3mm}
\psset{xunit=3mm}
\psset{yunit=3mm}
\begin{pspicture*}(40,16)(-25,0)
\multiput(-6.2,13.7)(2,-2){4}{$\bullet$}
\multiput(-8.2,11.7)(2,-2){6}{$\bullet$}
\multiput(-8.2,7.7)(2,-2){2}{$\bullet$}
\multiput(-10.2,9.7)(-2,-2){4}{$\bullet$}
\multiput(-10.2,5.7)(-2,-2){3}{$\bullet$}
\multiput(1.8,9.7)(4,0){5}{$\bullet$}
\multiput(1.8,5.7)(4,0){5}{$\bullet$}
\multiput(3.8,7.7)(4,0){5}{$\bullet$}
\multiput(3.8,3.7)(12,0){2}{$\bullet$}
\multiput(7.8,11.7)(4,0){2}{$\bullet$}
\multiput(21.8,9.7)(2,2){3}{$\bullet$}
\put(9.8,13.7){$\bullet$}
\put(-7,14.8){$\sigma_{16}$}
\put(9.45,14.8){$\s_8$}
\put(22.9,12.7){$H$}
\put(3.9,2.6){$\s_{11}''$}
\psline(-6,14)(0,8)
\psline(0,8)(2,10)
\psline(0,8)(2,6)
\psline(2,6)(10,14)
\psline(4,4)(12,12)
\psline(10,14)(12,12)
\psline(12,12)(18,6)
\psline(10,6)(14,10)
\psline(2,6)(4,4)
\psline(2,10)(4,8)\psline(4,8)(6,6)
\psline(6,10)(8,8)
\psline(8,12)(10,10)
\psline(10,10)(16,4)
\psline(8,8)(10,6)
\psline(16,4)(20,8)
\psline(14,6)(18,10)
\psline(18,10)(20,8)
\psline(20,8)(24,12)
\psline(24,12)(26,14)
\psline[linecolor=blue](2,2)(4,4)\psline[linecolor=blue](0,4)(2,6)
\psline[linecolor=blue](-2,6)(0,8)
\psline(-6,6)(-4,8)\psline[linecolor=blue](-4,8)(-2,10)
\psline(-14,2)(-6,10)\psline[linecolor=blue](-6,10)(-4,12)
\psline[linecolor=blue](-8,12)(-6,14)
\psline(-8,12)(2,2)
\psline(-10,10)(-6,6)\psline(-13,9)(-10,6)
\psline(-15,7)(-12,4)\psline(-17,5)(-14,2)
\psline(-16,4)(-8,12)
\put(1.6,.8){$q$}\put(-2.2,3){$qH$}
\end{pspicture*}
\end{center}

This diagram could be extended indefinitely on the left, with a period increasing degrees by $8$. Then it
would exactly encode the quantum Chevalley formula, in the sense that the product of any 
Schubert class by the hyperplane class $H$, would be the sum of the Schubert classes (possibly 
with some $q$-coefficient) connected to it  on its immediate left. We have drawn in blue the 
edges that introduce quantum corrections, so that the rightmost component of the diagram 
obtained by disconnecting the blue edges in just the classical Hasse diagram encoding the
classical Chevalley formula. In particular the diagram above, the {\it quantum Hasse diagram}, makes
it easy to compute powers of $H$ in the quantum cohomology ring. This will be useful a little later on. 

\smallskip
Over $\QQ[q]$, it was proved in \cite{cmp2} that the quantum cohomology ring is generated by the 
three classes $H$, $\sigma_8$ and $\sigma''_{11}$. The complex involution is thus determined by 
the images of these classes, which are sent by strange duality  to 
$$\iota(H)=12^{1/4}q^{-1} \sigma''_{11}, \quad\iota(\sigma_8)=q^{-2}\sigma_{16},
\quad \iota(\sigma''_{11})=12^{1/4}q^{-1}H.$$

\medskip\noindent{\it Proof of \ref{free}}. 
The first difficulty is to understand the $W(E_6)$-invariant polynomials on a 
Cartan subalgebra $\ft$ of $\fe_6$. For this, 
we consider the realisation of $\ft$ and $W(E_6)$ studied  in
\cite{mehta}. This is based on the observation that $\fe_6$ contains a 
full rank subalgebra isomorphic with $\fsl_6\times\fsl_2$. As modules 
over this subalgebra, $\fe_6$ and its minimal representation $J$, of dimension $27$, 
can be decomposed as 
\begin{eqnarray*}
\fe_6 &= &\fsl_6\times\fsl_2 \oplus \wedge^3\CC^6\otimes \CC^2, \\
J &= &(\CC^6)^*\otimes \CC^2\oplus \wedge^2\CC^6.
\end{eqnarray*}
We choose for our Cartan subalgebra $\ft$ of $\fe_6$ the product of Cartan subalgebras 
in $\fsl_6$ and $\fsl_2$. We thus have generators of $\ft^*$ given by $(x_i)_{1\leq i \leq 6}$ 
with $\sum x_i = 0$ (for the factor in $\fsl_6$), and $y$ (for the factor in $\fsl_2$). 
Moreover, the $W(E_6)$-invariants are generated by the obvious invariants given by the 
sums of powers of weights in $J$, that is,
$$I_k(x,y) = \sum_{1 \leq i , j \leq 6} (x_i + x_j)^k + \sum_i 
[(-x_i + y)^k + (-x_i - y)^k],$$
for $k \in \{2,5,6,8,9,12\}$. So, the equations we have to solve are
$I_k(x,y) = 0$ for $k \in \{2,5,6,8,9\}$ and $I_{12}(x,y) = c$, where $c$
is any non zero constant.

We claim that all the solutions of these equations, for a certain choice of $c$, 
 are given by the $W(E_6)$-orbit of the solution
$$
\left \{
\begin{array}{l}
x_2 = - x_1 = 1 \\
x_4 = - x_3 = i\\
x_6 = - x_5 = e^{i\pi/8} \\
y = e^{5i\pi/8},
\end{array}
\right .
$$
and that this $W(E_6)$-orbit is free. 
In fact, let us first check that this is indeed a solution. Since
$x_1 + x_2 = x_3 + x_4 = x_5 + x_6 = 0$, all $I_k$'s with odd $k$ vanish. 
For even $k$, $I_k$ reduces to a symmetric polynomial in 
$x_1^2,x_3^2,x_5^2$ and $y^2$. Therefore, to prove that $(x,y)$
annihilates $I_2,I_6$ and $I_8$, it is enough to show that
$$\sum_{i \in \{1,3,5,7\}} x_i^2 = 0, \quad
\sum_{i \in \{1,3,5,7\}} x_i^6 = 0, \quad 
\sum_{i \in \{1,3,5,7\}} x_i^8 = 0,$$ 
where we have set $x_7 := y$.
The first two equations hold because they are of odd degree in the $x_i^2$,
and the values of $x_i^2$ are two opposite pairs. To check that the last 
one also holds is a straighforward computation.

Let us now prove that the $W(E_6)$-orbit of our solution is free. To this end, let
$w \in W(E_6)$ such that $w(x,y) = (x,y)$. We know that $w$ can be represented
by a $7\times 7$ matrix $(a_{i,j})$  with rational coefficients. 
Since $x_1,x_3,x_5,y$ are independant over $\Q$, the diagonal coefficient $a_{7,7}$ 
of this matrix must be 1. Using the fact that
$w$ preserves $I_2$, we deduce that $w$ stabilises the vector
$(0,0,0,0,0,0,1)$. This means that $w$ belongs to the symmetric group $\cS_6\subset W(E_6)$,
acting by permutations of $x_1,\ldots ,x_6$.  But then  it
follows easily that $w$ is the identity.\qed

\bigskip\noindent
{\it Proof of \ref{sameinv}}. 
Let us first show that the complex conjugate of $\sigma_8$ is $\sigma_{16}$, as a warm-up. 
For this we recall that $\sigma_8$ is an idempotent of the quantum 
cohomology ring. More precisely, we have the following statement, which is a special
case of Corollary 5.2 in \cite{cmp2}. 

\begin{lemm}
The multiplication by $\sigma_8$ in $QH^*(E_6/P_1)$ sends any Schubert class to its 
translate by eight steps on the left in the quantum Hasse diagram. In particular
$\sigma_8^2 = \sigma_{16}$ and $\sigma_8^3 = q^2$.
\end{lemm}

In particular, we deduce that $\sigma_8$, considered as a function on 
$Z(E_6/P_1)$,
takes its values among the third roots 
of unity. But for such a root $\zeta$, we have $\overline\zeta=\zeta^2$, and therefore 
$$\overline{\sigma_8}=\sigma_8^2 = \sigma_{16}.$$

Let us now show that the complex conjugate of $H$ is $12^{1/4}\sigma''_{11}$.
This
is much more tricky, but the idea is the same as above:
suppose we can find a non zero polynomial $P$ such that $P(H)=0$; suppose that
we can also find a polynomial $Q$ such that $Q(\zeta)=\overline \zeta$ for 
any $\zeta$ root of $P$. 
Then we can conclude that $\overline H=Q(H)$. (That's exactly what we have   
done above with $P(\zeta)=\zeta^3-1$ and $Q(\zeta)=\zeta^2$). 

We first find an equation for $H$. For this we will find a dependance relation 
between $H^{25}$, $H^{13}$ and $H$ (recall that the index of the Cayley plane
is twelve).  We compute in the specialisation $QH^*(E_6/P_1)_{q=1}$. With the help of 
the quantum Chevalley formula, or equivalently, of the quantum Hasse diagram above, we get 
$$
\left \{
\begin{array}{l}
H^{13} = 78 \sigma_{13} + 57 H\\
H^{25} = 21060 \sigma_{13} + 15417 H,
\end{array}
\right .
$$
hence the equation $P(H) = H^{25} - 270 H^{13} - 27H = 0$. Note that 
the roots of $P$ are 
given by $\zeta=0$ or $\zeta^{12} = 135 \pm 78\sqrt 3 = (3 \pm 2 \sqrt 3)^3$.

\smallskip
The quantum Hasse diagram also gives us 
$$
\left \{
\begin{array}{l}
H^{11} = 33 \sigma'_{11} + 12 \sigma''_{11}\\
H^{23} = 8901 \sigma'_{11} + 3258 \sigma''_{11}
\end{array}
\right .
$$
Therefore, $234\sigma''_{11} = H^{11}(11H^{12} - 2967)$.

To check that $\overline H = 12^{1/4}\sigma''_{11}$, we therefore just need to
verify that for every root $\zeta$ of the polynomial $P(z)=z^{25} - 270 z^{13} - 27z$,
one has 
$234.12^{1/4}\overline \zeta = \zeta^{11}(11\zeta^{12} - 2967)$, 
or, equivalently,
$$234 \zeta\overline \zeta = \zeta^{12}(11\zeta^{12} - 2967).$$ 
This is clear if $\zeta=0$. If
$\zeta^{12} = 135 + 78 \sqrt 3$,  an explicit computation yields
$\zeta^{12}(11\zeta^{12} - 2967) = 702 +  234 \sqrt 3$, whereas
$12^{1/4}\zeta \overline \zeta   = \sqrt{(3+2\sqrt 3).2\sqrt 3} = 3 + \sqrt 3$.
The computation  is similar when $\zeta^{12} = 135 - 78 \sqrt 3$, and we can therefore
conclude that $\overline H = 12^{1/4} \sigma''_{11}$, as claimed. 

Of course we can deduce that, conversely, the complex conjugate 
of $\sigma''_{11}$ is $12^{-1/4}H$. Since $QH^*(E_6/P_1)$ is generated by 
$H,\sigma''_{11}$ and $\sigma_8$, this completes the proof. \qed

\subsection{The Freudenthal variety}

We recall some facts on the quantum cohomology ring of the Freudenthal variety  $E_7/P_7$. 
This is a $27$ dimensional variety, of index $18$. The quantum Hasse diagram is as follows: 

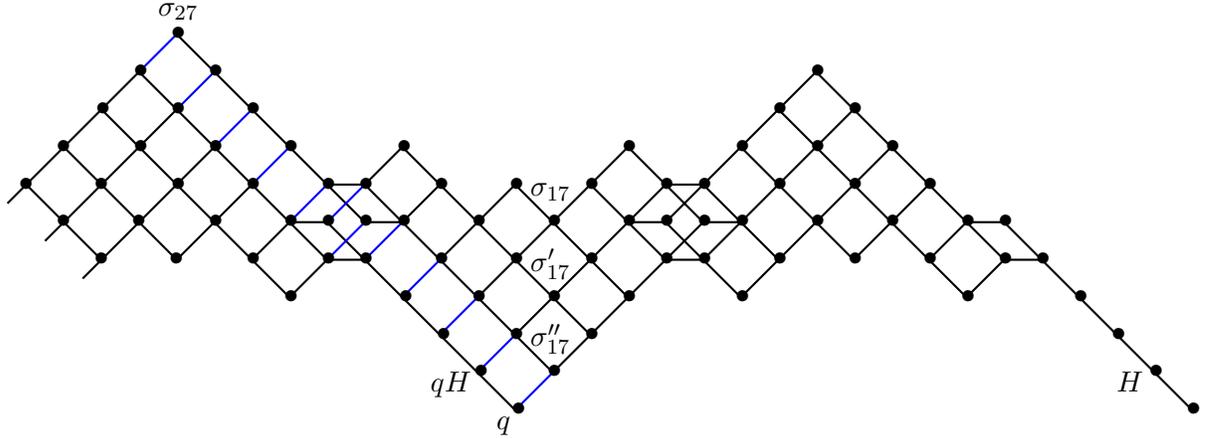
\begin{center}
\psset{unit=2.5mm}
\psset{xunit=2.5mm}
\psset{yunit=2.5mm}
\begin{pspicture*}(50,20)(-20,-6)
\psline(-8,18)(2,8)
\psline(-10,16)(0,6)\psline(2,6)(10,-2)
\psline(-12,14)(-2,4)\psline(-14,12)(-8,6)\psline(-16,10)(-12,6)
\psline(-10,16)(-17,9)\psline(-8,14)(-15,7)\psline(-6,12)(-13,5)
\psline(-4,10)(-8,6)\psline(-2,8)(-4,6)\psline(0,6)(-2,4)
\psline(2,8)(4,8)\psline(0,10)(2,10)\psline(0,6)(2,6)
\psline(2,10)(4,12)
\psline(2,10)(10,2)\psline(10,2)(12,0)
\psline[linecolor=blue](0,8)(2,10)\psline(0,8)(2,6)\psline(0,8)(-2,8)
\psline[linecolor=blue](2,6)(4,8)
\psline[linecolor=blue](4,4)(6,6)
\psline[linecolor=blue](6,2)(8,4)
\psline[linecolor=blue](8,0)(10,2)
\psline[linecolor=blue](10,-2)(12,0)
\psline[linecolor=blue](0,6)(2,8)
\psline[linecolor=blue](-2,8)(0,10)
\psline[linecolor=blue](-4,10)(-2,12)
\psline[linecolor=blue](-6,12)(-4,14)
\psline[linecolor=blue](-8,14)(-6,16)
\psline[linecolor=blue](-10,16)(-8,18)
\psline(4,8)(6,10)
\psline(6,6)(10,10)
\psline(4,12)(12,4)\psline(12,4)(14,2)
\psline(10,10)(14,6)\psline(14,6)(16,4)
\psline(16,12)(18,10)\psline(18,10)(20,8)
\psline(14,10)(16,8)\psline(16,8)(18,6)
\psline(8,4)(16,12)
\psline(10,2)(18,10)
\psline(12,0)(20,8)
\psline(16,8)(18,8)
\psline(18,6)(20,6)
\psline(18,10)(20,10)
\psline(20,8)(22,8)
\psline(18,8)(26,16)
\psline(18,8)(20,6)\psline(20,6)(22,4)
\psline(20,6)(28,14)
\psline(22,4)(30,12)
\psline(20,10)(22,8)\psline(22,8)(24,6)
\psline(22,12)(24,10)\psline(24,10)(28,6)
\psline(24,14)(26,12)\psline(26,12)(34,4)
\psline(26,16)(28,14)\psline(28,14)(36,6)
\psline(28,6)(32,10)
\psline(32,6)(34,8)\psline(34,8)(36,8)
\psline(34,4)(36,6)\psline(36,6)(38,6)
\psline(36,8)(46,-2)
\multiput(-8.3,17.7)(2,-2){6}{$\bullet$}
\multiput(-10.3,15.7)(2,-2){6}{$\bullet$}
\multiput(-12.3,13.7)(2,-2){6}{$\bullet$}
\multiput(-14.4,11.7)(2,-2){4}{$\bullet$}
\multiput(-16.4,9.7)(2,-2){3}{$\bullet$}
\multiput(3.8,3.7)(2,-2){4}{$\bullet$}
\multiput(1.7,9.7)(4,0){5}{$\bullet$}
\multiput(1.7,5.7)(4,0){5}{$\bullet$}
\multiput(3.7,7.7)(4,0){5}{$\bullet$}
\multiput(7.7,3.7)(4,0){3}{$\bullet$}
\multiput(3.7,11.7)(12,0){2}{$\bullet$}
\multiput(9.7,1.7)(4,0){2}{$\bullet$}
\put(11.7,-.3){$\bullet$}\put(-.3,7.7){$\bullet$}
\multiput(17.7,7.7)(4,0){5}{$\bullet$}
\multiput(19.7,9.7)(4,0){4}{$\bullet$}
\multiput(19.7,5.7)(4,0){5}{$\bullet$}
\multiput(21.7,11.7)(4,0){3}{$\bullet$}
\multiput(23.7,13.7)(4,0){2}{$\bullet$}
\multiput(25.7,15.7)(4,0){1}{$\bullet$}
\multiput(21.7,3.7)(12,0){2}{$\bullet$}
\multiput(35.7,7.7)(2,-2){6}{$\bullet$}
\put(42,-1){$H$}\put(-9,19){$\sigma_{27}$}
\put(10.8,9.4){$\sigma_{17}$}\put(10.8,5.4){$\sigma'_{17}$}\put(10.8,1.4){$\sigma''_{17}$}
\put(9,-3){$q$}\put(5.5,-1){$qH$}
\end{pspicture*}
\end{center}

Over $\QQ[q]$, the quantum cohomology ring of the 
Freudenthal variety is generated by 
the hyperplane class $H$ and the Schubert classes 
$\sigma_{17}$ and $\sigma_{27}$, 
the punctual class. The strange duality $\iota$ maps these classes to 
$$\iota(H)=3456^{1/9}q^{-1} \sigma_{17},\quad 
\iota(\sigma_{17})=3456^{-1/9}q^{-1}H, 
\quad \iota(\sigma_{27})=q^{-3}\sigma_{27}.$$

\medskip\noindent{\it Proof of \ref{free}}. 
In this case, we were not able to give explicitely any point in $Z(E_7)$, so
we give a more abstract argument. First recall that the 
$W(E_7)$-invariants polynomials on a Cartan subalgebra $\ft$ of $\fe_7$, 
are generated by homogeneous polynomials $J_l$ of degrees $l=2,6,8,10,12,18$.
To compute $Z(\fe_7)$, we will, as in the preceeding cases, give a
(non-explicit) example
of an element in this scheme, and show that it belongs to a free 
$W(E_7)$-orbit, so that $Z(\fe_7)$ is exactly this orbit and is reduced.

The orthogonal of $\omega_7$ in $\ft$ is a Cartan algebra for $\fe_6$; let
$I_k,k \in \{2,5,6,8,9,12\}$ denote $W(E_6)$-invariants of degree $k$ in 
$\omega_7^\bot$ which generate the algebra of invariants.
Let $u$ in $\ft$ be a point in $\omega_7^\bot$ such that
$I_2(u)=I_5(u)=I_6(u)=I_8(u)=0$ and $I_9(u) = 1$.

First we claim that $u \in Z(\fe_7)$. In fact, the 
$J_l,l\in\{2,6,8,10,12\}$ restrict on
$\omega_7^\bot$ to $W(E_6)$-invariants, therefore to polynomials in the
$I_k,k \in \{2,5,6,8,9,12\}$, and since 
$I_2(u)=I_5(u)=I_6(u)=I_8(u)=I_{12}(u)=0$, we must have $J_l(u) = 0$.

Let $w \in W(E_7)$ such that $w.u = u$.
Let $J$ denote the first fundamental representation of $\fe_6$,
and consider the characteristic polynomial of the 
action of $u$ on $J$~: since it is a $W(E_6)$-invariant of degree 27
and $I_2(u)=I_5(u)=I_6(u)=I_8(u)=I_{12}(u)=0$,
it must be of the form 
$$ Q(u,T) = \prod_{\eta} (T-\eta(u))=T^{27}+q_9T^{18}+q_{18}T^{9}+q_{27},$$
where the product is taken over the $27$ weights $\eta$ of $J$. 

Let $\eta_1$ be a weight of $J$ such that $x:=\eta_1(u) \not = 0$.
Let $\zeta$ be a primitive 9-th root of unity.
Since the set $\{\eta(u)\}$ is the set of roots of
$Q(u,T)$, which is a polynomial in $T^9$, 
there must exist weights $\eta_2,\ldots,\eta_6$
such that 
$$(\eta_1(u),\ldots,\eta_6(u))=
(x,x\zeta,x\zeta^2,\ldots,x\zeta^5).$$
Recall that the product $\prod_\theta (T-\theta)$
over the 6 primitive 9-th roots
of unity, a cyclotomic polynomial, is irreducible. Thus
$T^6 + T^3 + 1$ is the minimal polynomial of $\zeta$ and
the family $(1,\zeta,\ldots,\zeta^5)$ is free over $\Q$. Therefore,
$(x,x\zeta,\ldots,x\zeta^5)$ is also free
over $\Q$, and a fortiori $(\eta_1,\ldots,\eta_6)$ is free over $\Q$, and thus
$(\eta_i,\omega_7)$
form a base over $\Q$ of the weight vector space.

Thus setting 
$$w.\omega_7 = a \omega_7 + \sum_{1 \leq i \leq 6} a_i \eta_i,$$
we get
$0=\omega_7(w^{-1}.u)=\sum a_i \eta_i(u)$, so 
$\forall i \in \{1,\ldots,6\},a_i = 0$ and 
$w.\omega_7 = a \omega_7$.
By the same argument, we also have 
$w.\eta_i = \eta_i$ modulo $\omega_7$,
so $w$ preserves $\omega_7^\bot$ and is trivial on it.
Recall that in any Weyl group 
the only non trivial elements acting trivially on some hyperplane are 
the reflections $s_{\alpha}$,
where $\alpha$ is some root (see \cite{bou}, Chap. V, {\S} 3.2, 
Th{\'e}or{\`e}me 1 (iv)). 
But  $\omega_7$ is not a root, so the line it generates
is not preserved by any reflection. Hence we must have $w=id$. \qed

\medskip\noindent{\it Proof of \ref{sameinv}}. 
We know from Theorem 3.3 in \cite{cmp2} that $\sigma_{27}^2 = q^3$. When we specialize at $q=1$,  
we deduce that $\sigma_{27}$ is real, hence equal to its complex conjugate. 

Let us now show that the complex conjugate of $H$ is  $3456^{1/9}\sigma_{17}$. 
For this, we first express  $\sigma_{17}$
as a polynomial in $H$. Using the quantum Hasse diagram, we compute 
$$
\left \{
\begin{array}{l}
H^{17} = 78 \sigma_{17} + 442 \sigma'_{17} + 748 \sigma''_{17}\\
H^{35} = 2252088 \sigma_{17} + 12969160 \sigma'_{17} + 22121896 
\sigma''_{17}\\
H^{53} = 66396246672 \sigma_{17} + 382360744192 \sigma'_{17} +
652206892048 \sigma''_{17}.
\end{array}
\right .
$$
So we get the relation $\sigma_{17}=H^{17}Q(H^{18})$, where we have set  
$$Q(T) = \frac{4237743313}{721278} - \frac{33629825}{77976}T +
\frac{84371}{5770224} T^2.$$ 

In fact, these formulas can be proved simply using the periodicity of the
quantum Hasse diagram of $E_7/P_7$, and the fact that the restriction of the
multiplication by $H^{18}$ to the degree 8 part of
$QH^*(E_7/P_7)_{q=1}$ has matrix
$$
\left (
\begin{array}{ccc}
598  &  1710  &   1938\\
3420  &  9832  &  11172  \\
5814  &  16758  &  19066
\end{array}
\right )
$$
in the base $\sigma_{26},\sigma_8,\sigma'_8$. Moreover, this shows that if
$P(T) = 64 - 401808 T + 29496 T^2  - T^3$, the characteristic polynomial
of this matrix, then $H^8 P(H^{18}) = 0$.
Now, it is easy to check that  $P$ has three real roots, all positive. 
Moreover, our key point is the following 
\medskip

\noindent {\bf Fact.}\hspace*{1cm} 
$(T Q(T))^{18} = (3456 T)^2\quad \mathrm{modulo}\;\;P.$

\medskip
In fact, somewhat painful computations lead to
$(T Q(T))^2 \simeq -1/40071 \ T^2+799/1083\ T+11696/13357$,
$(T Q(T))^6 \simeq -40/13357\ T^2+34544/361\ T+8768/13357$, and finally
$(T Q(T))^{18} = (3456 T)^2\quad \mathrm{mod}\;\;P.$
Therefore, if $\zeta$ is such that $P(\zeta^{18}) = 0$, then
$(\zeta^{18}Q(\zeta^{18}))^{18} = (3456 \zeta^{18})^2$, and since
$\zeta^{18}$ is real, this is equal to $3456^2  (\zeta \overline \zeta)^{18}$. 
We have noticed that $\zeta^{18}$, being  a root of $P$, is positive. One also 
checks that $Q(\zeta^{18}) \ge 0$. So we can extract $18$-th roots, 
and deduce that 
$$\overline \zeta = 3456^{1/9}\zeta^{17}Q(\zeta^{18}).$$
Note that this also holds for $\zeta=0$, so 
this is exactly what we need to conclude that the complex conjugate of $H$ is
$$\overline H = 3456^{1/9}H^{17}Q(H^{18}) 
= 3456^{1/9}\sigma_{17}.$$
Of course we can deduce that the complex conjugate of $\sigma_{17}$ is the 
expected multiple of $H$, and this concludes the proof. \qed 

\section{A non reduced example}


In this section we consider the symplectic Grassmannian $G_{\omega}(2,6)$. Its 
dimension is $7$ and its index is $5$. Its Schubert classes are indexed by pairs
$(a|\alpha)$, where $0\le a\le 2$ and $\alpha$ is a strict partition whose parts
are not bigger than $3$, and whose length is at most $a$. The degree of the 
corresponding class is $a$ plus the sum of the parts of $\alpha$. 
The Hasse diagram is the following: 
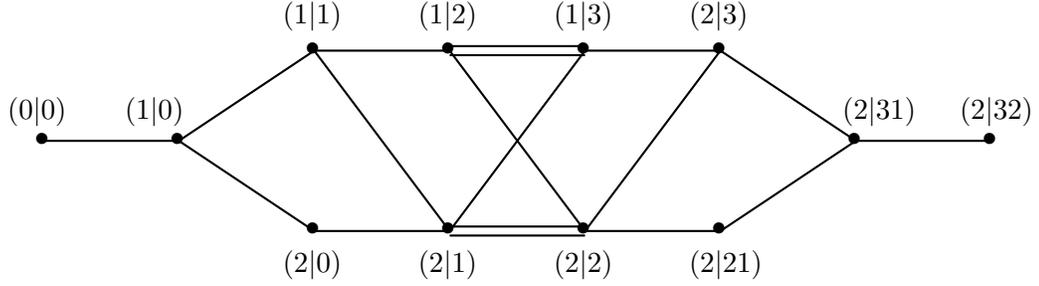
\begin{center}
\psset{unit=6mm}
\psset{xunit=6mm}
\psset{yunit=6mm}
\begin{pspicture*}(30,8)(-6,0)
\multiput(-2.2,3.9)(3,0){2}{$\bullet$}
\multiput(3.8,5.9)(3,0){4}{$\bullet$}
\multiput(3.8,1.9)(3,0){4}{$\bullet$}
\multiput(15.8,3.9)(3,0){2}{$\bullet$}

\psline(-2,4)(1,4)
\psline(1,4)(4,2)
\psline(1,4)(4,6)
\psline(4,2)(7,2)
\psline(4,6)(7,6)
\psline(4,6)(7,2)
\psline(7,6)(10,2)
\psline(7,2)(10,6)
\psline(7,5.9)(10,5.9)
\psline(7,6.1)(10,6.1)
\psline(7,1.9)(10,1.9)
\psline(7,2.1)(10,2.1)
\psline(10,2)(13,2)
\psline(10,2)(13,6)
\psline(10,6)(13,6)
\psline(13,2)(16,4)
\psline(13,6)(16,4)
\psline(16,4)(19,4)

\put(-2.8,4.5){$(0|0)$}
\put(-.2,4.5){$(1|0)$}
\put(3.3,6.6){$(1|1)$}
\put(3.3,1.1){$(2|0)$}
\put(6.3,6.6){$(1|2)$}
\put(6.3,1.1){$(2|1)$}
\put(9.3,6.6){$(1|3)$}
\put(9.3,1.1){$(2|2)$}
\put(12.3,6.6){$(2|3)$}
\put(12.3,1.1){$(2|21)$}
\put(15.7,4.5){$(2|31)$}
\put(18.3,4.5){$(2|32)$}

\end{pspicture*}
\end{center}

Tamvakis proved in \cite{T} that the quantum cohomology ring of $G_{\omega}(2,6)$
is generated by the special Schubert classes $\sigma_k=\sigma_{(1|k-1)}$, for $1\le k\le 4$, 
with the following relations:
$$\s_2^2-2\s_1\s_3+2\s_4=\s_3^2-2\s_2\s_4+q\s_1=0,$$

$$\begin{pmatrix} \s_1 &\s_2 &\s_3 \\ 1 & \s_1& \s_2 \\ 0& 1 & \s_1 \end{pmatrix}=0, 
\qquad
\begin{pmatrix} \s_1 &\s_2 &\s_3 &\s_4\\ 1 & \s_1& \s_2 &\s_3 \\ 
0& 1 & \s_1 &\s_2 \\ 0 & 0 & 1 & \s_1\end{pmatrix}=0. $$
The last two relations allow to express $\s_3$ and $\s_4$ as polynomials in $\s_1,\s_2$. 
Plugging these expressions in the first two relations, and letting $q=1$, we get 
$$\s_2(3\s_2-2\s_1^2) = 2\s_2^2\s_1^2-2\s_2\s_1^4+\s_1^6-2\s_2^3+\s_1=0.$$
These two equations define the finite scheme $Z(G_{\omega}(2,6))$. This scheme has length 12, 
and is the union 
of a subscheme $Z_1$, the union of the 5 simple points given by $\s_2=0, \s_1^5=-1$, 
a subscheme $Z_2$, the union of the other 5 simple points given 
by $3\s_2=2\s_1^2, \s_1^5=27$, and a length two scheme $Z_0$ supported at the origin, 
with Zariski tangent space $\s_1=0$. The quantum cohomology ring of $G_{\omega}(2,6)$
at $q=1$, identifies with the algebra of functions on that scheme $Z$. It is convenient 
to identify each Schubert class with such a function, given by a triple of functions 
on the subschemes $Z_0, Z_1, Z_2$. For this we just need the quantum Chevalley formula, 
and the result is as follows:

$$\begin{array}{cccccccccccc}
    & \s_{(1|0)} &\s_{(1|1)}      &\s_{(0|2)}     &\s_{(2|1)}     &\s_{(1|2)}     & \s_{(2|2)} &
\s_{(1|3)}   &\s_{(2|3)} &\s_{(2|21)}&\s_{(2|31)} &\s_{(2|32)}     \\
 & & & & & & & & & & &\\
Z_0 & 0          &\epsilon        & -\epsilon     & 0             & 0             & 0 &
0            & -1        & 1         & 0          & -\epsilon       \\
Z_1 & s          & 0              & s^2           & s^3           & -s^3          & s^4 &
 -s^4         & 0         & -1        & -s         & -s^2           \\
Z_2 & s          & \frac{2}{3}s^2 & \frac{1}{3}s^2& \frac{1}{3}s^3& \frac{1}{3}s^3& \frac{1}{9}s^4 &
\frac{1}{9}s^4& 2         & 1         & s          & \frac{1}{3}s^2 
\end{array}$$
 
\bigskip 
Note in particular that $\s_{(2|32)}-\s_{(1|1)}+\s_{(0|2)}$ is nilpotent: it generates the
radical of the quantum cohomology ring at $q=1$. In fact, even before specializing $q$, we 
have 
$$(\s_{(2|32)}-q\s_{(1|1)}+q\s_{(0|2)})^2=0\qquad\mathrm{in}\quad QA^*(G_{\omega}(2,6)).$$
It would  now be completely straightforward to write down the multiplication table of Schubert classes 
in the quantum cohomology ring. But the point we want to stress is the following: if we consider a degree
two class, $\s_{(1|1)}$ or $\s_{(0|2)}$, its complex conjugate is non zero on $Z_0$, and therefore it 
cannot be expressed as a linear combination of the degree three  classes $\s_{(2|1)}$ and $\s_{(1|2)}$,
which vanish on $Z_0$. It is therefore impossible to lift the complex conjugation to the localized quantum
cohomology ring of $G_{\omega}(2,6)$, in such a way that a degree $k$  class is mapped to a class of
degree $-k$. The existence of a non reduced point, although as simple as possible, definitely prevents
us from doing that. 

The fact that the first non (co)minuscule example leads to a non semi-simple quantum cohomology 
algebra suggests the following 

\medskip
 {\bf Conjecture}. {\it Consider  a rational homogeneous variety $G/P$ with Picard number one. 
Then $QA^*(G/P)_{q=1}$, its quantum cohomology algebra specialized  at $q=1$, is semi-simple, if and 
only if $G/P$ is minuscule or cominuscule.} 
 
\section{Vafa-Intrilagator type formulas}

\subsection{The quantum Euler class}

Abrams introduced in  \cite{ab}, for any projective variety $X$ (in fact, in a more general setting)  
a {\it quantum  Euler class} $e(X)$.
Let us denote by $pt$ the element of $W^P$ which defines the punctual 
class $\s(pt)$.
Let $\varphi : QH^*(X,\C) \simeq H^*(X,\C) \otimes_\C \C[q] \rightarrow \C[q]$
be defined by $\varphi(\sum P_\lambda(q)\sigma(\lambda)) = P_{pt}$. Then
the bilinear map $(x,y) \mapsto \varphi(x * y)$ defines a Frobenius algebra
structure on $QH^*(X,\C)$ over $\C[q]$. Let $(e_i)$ be a base of $QH^*(X,\C)$
over $\C[q]$ and $(e_i^*)$ the dual base. Our bilinear form identifies
$QH^*(X,\C)$ and its dual; let $e_i^\sharp \in QH^*(X,\C)$ correspond to $e_i^*$.
The Euler class is defined by $e(X) := \sum e_i.e_i^\sharp$.
By \cite{ab}, Theorem 3.4), in a given specialization of the quantum 
parameters, the
invertibility of $e(X)$ is equivalent to the semi-simplicity of that algebra.
We will identify that class explicitely when $X=G/P$ is any homogeneous variety.
 
Recall that the invariant algebras 
$$\QQ[\ft]^{W_P}=\QQ[X_1,\ldots , X_r] \quad {\rm and} \quad  
 \QQ[\ft]^W=\QQ[Y_1,\ldots , Y_r]$$
are polynomial algebras over certain homogeneous invariants, by the famous Chevalley theorem. 
The classical and quantum cohomology rings of $G/P$ are 
$$H^*(G/P)_\QQ=\QQ[\ft]^{W_P}/(Y_1,\ldots ,Y_r)\quad\mathrm{and}\quad 
QA^*(G/P)_\QQ=\QQ[\ft]^{W_P}[q]/(Y_1(q),\ldots ,Y_r(q)),$$
where $Y_1(q),\ldots ,Y_r(q)$ are certain homogeneous deformations of the classical 
relations $Y_1,\ldots ,Y_r$. 


Denote by $\Phi(G/P)=\Phi^+-\Phi^+_P$ the set of positive roots of $\fg$ that are not roots of $\fp=Lie(P)$. 
This is also the set of weights in $\fg/\fp$, the $P$-module whose associated vector bundle on $G/P$ is the 
tangent bundle. Each root in $\Phi(G/P)$ defines a linear functional on $\ft$, and since $\Phi(G/P)$ is 
preserved by $W_P$, the product of the roots it contained defines a $W_P$-invariant polynomial function 
on $\ft$, hence a class in the cohomology ring $QA^*(G/P)_{q=1}$. 

\begin{prop}
The quantum Euler class of $G/P$ is 
$$e(G/P)=\prod_{\alpha\in\Phi(G/P)}\alpha.$$
\end{prop}

\proof 
Chose a basis $t_1,\ldots ,t_r$ of $\ft^*$. There is a constant $c\ne 0$ such that 
$$\det\Big(\frac{\partial Y_p}{\partial t_q}\Big)_{1\le p,q\le r} = c\prod_{\alpha\in\Phi^+}\alpha.$$
Indeed both terms are proportional to the minimal degree anti-invariant of the Weyl group
(see \cite{bou}, Chap. V, {\S} 5.4, Proposition 5).   Similarly, for $W_P$, we get
$$\det\Big(\frac{\partial X_p}{\partial t_q}\Big)_{1\le p,q\le r} = c_P\prod_{\alpha\in\Phi^+_P}\alpha,$$
and both terms are proportional to the minimal degree anti-invariant of $W_P$.  
We deduce that 
$$J:=\det\Big(\frac{\partial Y_p}{\partial X_q}\Big)_{1\le p,q\le r} = d_P\prod_{\alpha\in\Phi(G/P)}\alpha,$$
for some non zero constant $d_P$. 

Now, it was proved in \cite{ab}, Proposition 6.3, that in our situation, $J$ and the quantum Euler class 
$e(G/P)$ coincide up to some invertible class
in $QH^*(G/P)$. Moreover, since $J$ is homogeneous 
of degree $\#\Phi(G/P)=\dim (G/P)$, this invertible class is in fact a constant (see the 
proof of Proposition 6.4). So $e(G/P)$ must be a constant multiple of $\prod_{\alpha\in\Phi(G/P)}\alpha$. 

We claim that this constant is one. Indeed, the quantum Euler 
class is a deformation of the classical Euler class (\cite{ab}, see the claim 
after Proposition 5.1). 
But the latter coincides with $c_{top}(T_{G/P})$, the top Chern class of the tangent space, 
which in the classical setting is precisely given by the product of the roots in $\Phi(G/P)$. 
This implies our claim and the proof is complete. \qed

\subsection{The Gromov-Witten invariants}

For a class $\sigma$ in $QH^*(G/P)_{q=1}$, 
denote by $\trace(\sigma)$ the trace of the multiplication operator by $\sigma$. In terms 
of our finite set $Z(G/P)$, we have 
$$\trace(\sigma)=\sum_{z\in Z(G/P)}\sigma(z).$$


Recall that we denote by $pt$ the element of $W^P$ which defines the punctual 
class $\s(pt)$. Recall that the Schubert class which is Poincar{\'e} dual to $\s(\mu)$ is 
denoted $\s(p(\mu))$. 

\begin{lemm}
For any Schubert classes $\sigma(\lambda)$ and $\sigma(\mu)$ on $G/P$, we have
$$ \trace\Big(\frac{\sigma(\lambda)}{e(G/P)}\Big) = \delta_{\lambda,pt}, \qquad
 \trace\Big(\frac{\sigma(\lambda)\sigma(\mu)}{e(G/P)}\Big) =  \delta_{\lambda,p(\mu)}.$$
\end{lemm}
 
\proof 
The first equality is a general fact in Frobenius algebras.
Recall that for $x \in QH^*(X,\C)$ we denoted $\varphi(x)$ the coefficient of
$x$ in $\sigma_{pt}$. We want to show that
$\trace(x) = \varphi(e(G/P)x)$ for $x \in QH^*(G/P,\C)$. But we have
$\varphi(e(G/P)x) = \sum_i \varphi(e_i^\sharp e_i x)
= \sum e_i^*(e_i x) = \trace(x)$ (the second equality follows from the
definition of $e_i^\sharp$).
\smallskip

We deduce that $\trace(\sigma(\lambda)\sigma(\mu)/e(G/P))$ is the coefficient of $\sigma(pt)$ in the 
quantum product $\sigma(\lambda)*\sigma(\mu)$, that is, the Gromov-Witten invariant $I(\sigma(\lambda),
\sigma(\mu),1)$. By the associativity of the quantum product, this is also the coefficient of 
the Poincar{\'e} dual class $\sigma(p(\mu))$ inside $\sigma(\lambda)*1=\sigma(\lambda)$. \qed

\medskip
Our Vafa-Intriligator type formula follows immediately (up to the 
identification of $J=e(G/P)$, this formula appears in \cite[Section 4]{st})~:

\begin{prop}
The three points genus zero Gromov-Witten invariants of $G/P$ can be computed as 
$$I(\sigma(\lambda),\sigma(\mu),\sigma(\nu))=\trace\Big(\frac{\sigma(\lambda)\sigma(\mu)\sigma(\nu)}{e(G/P)}\Big)
=\sum_{z\in Z(G/P)}\frac{\sigma(\lambda)(z)\sigma(\mu)(z)\sigma(\nu)(z)}{
\prod_{\alpha\in\Phi(G/P)}\alpha(z)}.$$
\end{prop}

Note that we have not indicated the degree $d$ of this invariant. In fact it is determined by the  
Schubert classes involved, through the relation $\deg\s(\lambda)+\deg\sigma(\mu)+\deg \sigma(\nu)
=\dim(G/P)+\mathrm{ind}(G/P)d$.

\smallskip
In fact, following \cite{st} we have a general formula for the genus $g$  Gromov-Witten invariants:
$$I^g(\sigma(\lambda),\sigma(\mu),\sigma(\nu))=\trace\Big(e(G/P)^{g-1}\sigma(\lambda)\sigma(\mu)\sigma(\nu)\Big).$$

It would be interesting to understand better the quantum Euler class of a (co)minuscule $G/P$. 
We conjecture that, as functions on $Z(G/P)$, we should have
$$e(G/P)=|e(G/P)|\sigma(pt).$$
Recall that the punctual class $\sigma(pt)$ is always an idempotent, so that it takes its  values 
among the roots of unity. So the punctual class would give the argument of the quantum Euler class. 
We have checked this on Grassmannians and quadrics, but we have 
no explanation of this intriguing relation.

\bigskip\noindent
Pierre-Emmanuel {\sc Chaput}, Laboratoire de Math{\'e}matiques Jean Leray, 
UMR 6629 du CNRS, UFR Sciences et Techniques,  2 rue de la Houssini{\`e}re, BP
92208, 44322 Nantes cedex 03, France. 

\noindent {\it email}: pierre-emmanuel.chaput@math.univ-nantes.fr

\medskip\noindent
Laurent {\sc Manivel}, 
Institut Fourier, UMR 5582 du CNRS,  Universit{\'e} de Grenoble I, 
BP 74, 38402 Saint-Martin d'H{\'e}res, France.
 
\noindent {\it email}: Laurent.Manivel@ujf-grenoble.fr

\medskip\noindent
Nicolas {\sc Perrin},  Institut de Math{\'e}matiques,  
Universit{\'e} Pierre et Marie Curie, Case 247, 4 place Jussieu,  
75252 PARIS Cedex 05, France.

\noindent {\it email}: nperrin@math.jussieu.fr
\end{document}